\documentclass[reqno,11pt]{amsart}

\usepackage[dvips,pdftex]{epsfig}
\usepackage{amssymb}
\usepackage{amsmath}
\usepackage{amscd}
\usepackage{amsthm}
\usepackage[mathscr]{eucal}
\usepackage[all]{xy}

\newtheorem{definition}{Definition}
\newtheorem{theorem}{Theorem}
\newtheorem{lemma}{Lemma}
\newtheorem{corollary}{Corollary}

\DeclareMathOperator{\tr}{Tr}

\DeclareMathOperator{\E}{E}
\DeclareMathOperator{\LL}{L}
\DeclareMathOperator{\HH}{H}
\DeclareMathOperator{\K}{K}
\DeclareMathOperator{\TT}{T}
\DeclareMathOperator{\res}{Res}

\newcommand{\ff}{\boldsymbol{\zeta}_{\gamma,\alpha}}
\newcommand{\ml}{\mathsf{L}}
\newcommand{\mm}{\mathsf{M}}
\DeclareMathOperator{\End}{End}
\newcommand{\hap}{\widehat{p}}

\DeclareMathOperator{\pic}{Pic}
\DeclareMathOperator{\spa}{span}

\newcommand{\mq}{\mathsf{Q}}

\begin{document}

\title[The linear flows in the space of Krichever-Lax matrices]{The linear flows in the space of Krichever-Lax matrices over an algebraic curve}
\author[Taejung Kim]{Taejung Kim}

\address{Korea Institute for Advanced Study\\
207-43 Cheongyangri-dong\\
Seoul 130-722, Korea}

\thanks{The author would like to express his sincere gratefulness to Prof.~William Goldman, Prof.~Serguei Novikov, and Prof.~Niranjan Ramachandran for giving him their insights and valuable suggestions while preparing this manuscript.}

\keywords{Krichever-Lax matrix,  Lax representation,  Hitchin system,  Eigenvector mapping, Isospectral deformation, Euler sequence, Elementary transformation, Sheaf cohomology, Linearizing flow.}

\subjclass[2000]{14F05, 14H70,  32L10, 58F07}
\date{\today}

\email{tjkim@kias.re.kr}

\begin{abstract}
In \cite{kri02}, I.~M. Krichever invented the space of matrices parametrizing the cotangent bundle of moduli space of stable vector bundles over a compact Riemann surface, which is named as the Hitchin system after the investigation \cite{hit87}. We study a necessary and sufficient condition for the linearity of flows on the space of Krichever-Lax matrices in a Lax representation in terms of cohomological classes using the similar technique and analysis from the work \cite{grif85} by P.~A. Griffiths.
\end{abstract}

\maketitle

\section{Introduction}

In N. Hitchin's investigation \cite{hit87}, the dynamics of Hamiltonians on the cotangent bundle of the moduli space of stable vector bundles on a compact Riemann surface is characterized by straight line flows. Indeed, it is a basic distinction between algebraically completely integrable systems and completely integrable systems. The essence of this characterization in \cite{hit87} comes from the existence of a larger symplectic manifold containing the cotangent bundle where each fiber, an open set of the Jacobi variety of a spectral curve, is naturally compactified. The extension of Hamiltonian vector fields to the larger symplectic manifold is equivalent to the straightness of the associated Hamiltonian flows, since each fiber is a complex torus.

In the space of Krichever-Lax matrices, a priori there is no symplectic structure nor Poisson structure defined. Due to this reason, we do not have any Hamiltonian dynamics yet. The starting point of \cite{kri02} is to define the dynamics of system on the space of Krichever-Lax matrices in terms of what is called a Lax representation:
\begin{equation}\label{de5lax}
\frac{d}{dt}\ml_t=[\mm_t,\ml_t].
\end{equation}
Note that $\mm$ is a function of $\ml$. The matrix $\mm$ characterizes the dynamics of flows in the space of Krichever-Lax matrices. I.~M. Krichever gives the condition on $\mm$ when the flows of the Lax representation become straight (\emph{Theorem 2.1, Theorem 2.2} in \cite{kri02}). Moreover, he constructs a symplectic structure on the space of Krichever-Lax matrices 
and shows that the straight line flows coming from the Lax representation indeed are Hamiltonian flows. That is, they define Hamiltonians associated with the symplectic structure.

In the meantime, the dynamics of Lax representation~\eqref{de5lax} is completely described by $\mm_t$ up to addition of a polynomial $\mathsf{P}(\ml_t)$ or a commuting element $\mq_t$ with $\ml_t$. It seems natural that this ambiguity in a Lax representation can be well encoded in cohomology classes. Thus, we will characterize the straightness of flows in terms of cohomology classes. Note that a priori the flows in a Lax representation are not necessarily straight line flows. In order to describe the Hitchin system using a Lax representation, we need a special condition on $\mm$. In \cite{grif85}, P.~A. Griffiths gave a necessary and sufficient condition where the flows from a Lax representation are straight in the case of spectral curves over $\mathbb{P}^1$. A similar question for the Hitchin system has not been answered yet in the author's knowledge (compare this with \cite{gat92}). In this paper, we will investigate this question and give an answer which will also verify that the choice of $\mm$ in \cite{kri02} indeed induces a straight flow. The main results are stated in Theorem~\ref{thco1}, Theorem~\ref{cothe1}, and Corollary~\ref{cohomco3}.

Let us summarize our main results. The flows associated with~\eqref{de5lax} can be seen as isospectral deformations of line bundles $\LL_t$. That is, they stay in the Jacobi varieties of isospectral objects, commonly named spectral curves $\widehat{\mathfrak{R}}$, of a compact Riemann surface $\mathfrak{R}$. In Theorem~\ref{cothe1}, 
we characterize how the tangent vector of a flow induces a class $\rho(\mm_t)$, so-called a \emph{residue section} (see Definition~\ref{rede2}), in a sheaf cohomology group associated with a skyscraper sheaf which is gotten by a lifting divisor of a multiple of the canonical divisor $K$ of $\mathfrak{R}$ as follows:
$$\frac{d}{dt}\LL_t=\partial\rho(\mm_t) \in\HH^1(\widehat{\mathfrak{R}},\mathcal{O}_{\widehat{\mathfrak{R}}}).$$
The linearity of tangent flows is equivalent to that the time-derivative of residue section is zero up to cohomology groups induced by exact sequence~\eqref{exsk1} of elementary transformation. This is indeed a generalization of the Griffiths' result \cite{grif85}, which says that the linearity of flows on the Jacobi variety of a spectral curve over $\mathbb{P}^1$ is characterized by the time-derivative of a cohomology class associated with a skyscraper sheaf of lifting divisor of the infinity and zero of $\mathbb{P}^1$. In Section~\ref{secch4}, we also explain that the criteria in \cite{kri02} on the linearity of Hamiltonian vector fields expounded by Krichever nicely fit into the scheme of our investigation.

\section{The Hitchin systems and the spaces of Krichever-Lax matrices}

We will briefly give basic definitions in Hamiltonian dynamics. For more detail, we refer to \cite{can01}. Let $M$ be a symplectic manifold with a symplectic form $\omega$. A \emph{Hamiltonian vector field} $X_H$ associated to the symplectic form $\omega$ and a smooth function $H$ on $M$ is defined by
$$dH=\iota(X_H)\omega.$$
We will call $H$ a \emph{Hamiltonian or Hamiltonian function}. The \emph{Poisson bracket} $\{\,,\,\}$ (\emph{p.108} in \cite{can01}) associated to the symplectic form is defined by
$$\{H,G\}=X_HG.$$
Two functions $H,G$ are said to be \emph{Poisson commutative} if
$$\{H,G\}=0.$$
Note that the maximal number of linearly independent Hamiltonians on a symplectic manifold $M$ of dimension $2n$ is $n$. Accordingly, we say that a symplectic manifold $M$ of dimension $2n$ is a \emph{completely integrable system} if it has $n$ linearly independent Hamiltonians $H_1,\dots,H_n$ generically, i.e., 
$$dH_1\wedge\cdots\wedge dH_n\ne0\text{ generically}.$$ 
If $M$ is a completely integrable system, we may define a map
$$\mathbf{H}:M^{2n}\to\mathbb{C}^n\text{ by }\mathbf{H}(m)=(H_1(m),\dots,H_n(m)).$$
This is a special case of a momentum map (\emph{p.133} in \cite{can01}) in symplectic geometry. 
Indeed, it is a momentum map for the action of an abelian group, i.e., a complex torus. The primary dynamical system to study in this paper is presented as follows:

\begin{definition}(\emph{p.96} in \cite{hit87})
A dynamical system is said to be an algebraically completely integrable system if 

\begin{itemize}
\item[1]
it is a completely integrable system
\item[2]
a generic fiber of $\mathbf{H}$ is an (Zariski) open set of an abelian variety
\item[3]
each Hamiltonian flow of $X_{H_i}$ is linear on a generic fiber.
\end{itemize}
\end{definition}

In \cite{hit87}, Hitchin proves that for a moduli space $\mathscr{N}$ of stable holomorphic vector bundles of rank $l$ over a compact Riemann surface $\mathfrak{R}$ of genus greater than $1$, the cotangent bundle $\TT^\ast\mathscr{N}$ is an algebraically completely integrable system, so-called a Hitchin system (see \cite{don96,hit87,kim07} for details). The main part of the proof builds up on an observation that 
a generic fiber of the Hitchin map $\mathbf{H}:\TT^\ast\mathscr{N}\to\bigoplus^{k}_{i=1}\HH^0(\mathfrak{R},\K^{d_i}_{\mathfrak{R}})$ defined by 
invariant polynomials
$$\mathbf{H}(\Phi_{[A]}(z))=\Big(h_1(\Phi_{[A]}(z)),\dots,h_k(\Phi_{[A]}(z))\Big)$$
is an open set of the Jacobi variety of a spectral curve $\widehat{\mathfrak{R}}$ associated with a Higgs field $\Phi_{[A]}$ at $[A]\in\mathscr{N}$ where
$$\widehat{\mathfrak{R}}=\{\lambda_z\in \K_\mathfrak{R}\mid \pi^\ast\det(\lambda_{z}\cdot I_{l\times l}-\Phi_{[A]}(z))=0\}.$$
Here $\K_\mathfrak{R}$ is the canonical bundle of $\mathfrak{R}$ and the invariant polynomials $h_1,\dots,h_k$ are the
coefficients of the characteristic polynomial of a Higgs field.

In the meantime, an explicit parametrization of the cotangent bundle of a moduli space of stable vector bundles over a compact Riemann surface is investigated in \cite{kri02}: Let $\E_{\gamma,\alpha}$ be a holomorphic vector bundle of rank $l$ on a compact Riemann surface $\mathfrak{R}$ of genus $g>1$ associated with Tyurin parameters $(\gamma,\alpha)=\Big\{\gamma_j,\boldsymbol{\alpha}_j\Big\}_{j=1}^{lg}\in
\mathcal{S}^{lg}(\mathfrak{R}\times\mathbb{P}^{l-1})$ (see \cite{kim07} or \cite{tyu65,tyu66} for details). A global section $\ff(z)\in\mathcal{F}_{\gamma,\alpha}$ can be written as a vector-valued meromorphic function on $\mathfrak{R}$: Let $\boldsymbol{\alpha}_j=(\alpha_{1,j},\dots,\alpha_{l-1,j},1)\in\mathbb{C}^l$ for $j=1,\dots,lg$. Then a local expression is given by

\begin{equation}\label{kriform1}
\ff(z)=\frac{c_j\boldsymbol{\alpha}_j}{z-z(\gamma_j)}+O(1)\text{ where }c_j\in\mathbb{C}.
\end{equation}
From the Riemann-Roch theorem and given constraint~\eqref{kriform1}, we have
$$\dim_\mathbb{C}\HH^0(\mathfrak{R},\mathcal{F}_{\gamma,\alpha})\geq l(lg-g+1)-lg(l-1)=l.$$
Such vector bundles $\E_{\gamma,\alpha}$ with mutually distinct $\gamma_j$ for $j=1,\dots,lg$ and satisfying $\dim_\mathbb{C}\HH^0(\mathfrak{R},\E_{\gamma,\alpha})=l$ form an open set $\mathcal{M}'_0$ of $\mathcal{S}^{lg}(\mathfrak{R}\times\mathbb{P}^{l-1})$. 

\begin{definition}\label{laxde1}
A Krichever-Lax matrix associated to Tyurin parameters $(\gamma,\alpha)$ and a canonical divisor $K$ of a compact Riemann surface $\mathfrak{R}$ of 
genus $g>1$ is a matrix-valued meromorphic function $\ml(p;\gamma,\alpha)$ with at most simple poles at $\gamma_i$ and poles at $K$ satisfying the following conditions: There exist $\boldsymbol{\beta}_j\in\mathbb{C}^{l}$ and $\kappa_j\in\mathbb{C}$ for $j=1,\dots,lg$ such that a local expression in a neighborhood of $\gamma_j$ is given by
$$\ml(p;\gamma,\alpha)=\frac{\ml_{j,-1}(\gamma,\alpha)}{z(p)-z(\gamma_j)}+\ml_{j,0}(\gamma,\alpha)+O\big((z(p)-z(\gamma_j))\big)$$
with the following two constraints
\begin{itemize}
\item[1.]
$\ml_{j,-1}(\gamma,\alpha)=\boldsymbol{\beta}_j^{T}\cdot\boldsymbol{\alpha}_j$, i.e., of rank $1$  and it is traceless
$$\tr\ml_{j,-1}=\boldsymbol{\alpha}_j\cdot\boldsymbol{\beta}_j^{T}=0.$$
\item[2.]
$\boldsymbol{\alpha}_j$ is a left eigenvector of $\ml_{j,0}$
$$\boldsymbol{\alpha}_j\ml_{j,0}(\gamma,\alpha)=\kappa_j\boldsymbol{\alpha}_{j}.$$
\end{itemize}
Let us denote the set of Krichever-Lax matrices associated to Tyurin parameters $(\gamma,\alpha)$ and a canonical divisor $K$ by $\mathcal{L}_{\gamma,\alpha}^{K}$. Note that we will also call a Krichever-Lax matrix a Lax matrix following the terminology in \cite{kri02} for simplicity.
\end{definition}

The two constraints imply that a Lax matrix can be thought as a \textit{Higgs field} $\ml(p;\gamma,\alpha)\otimes\omega$, i.e., a global section of $\End(\E_{\gamma,\alpha})\otimes \K$: In a neighborhood of $\gamma_j$, 
the first and the second condition respectively imply

\begin{equation}\label{krieq2}
\begin{aligned}
\ff(z)\ml_{j,-1}(\gamma,\alpha)&=O(1)\\
\ff(z)\ml_{j,0}(\gamma,\alpha)&=\kappa_j\ff(z).
\end{aligned}
\end{equation}
Since we are assuming the divisor $K$ of $\omega$ does not intersect with $\{\gamma_j\}_{j=1}^{lg}$, 
we may conclude that $\ml(p;\gamma,\alpha)\otimes\omega$ is a global section of $\End(\E_{\gamma,\alpha})\otimes \K$. 
The dimension of the space of Lax matrices is
$$\dim_\mathbb{C}\mathcal{L}^K=l^2(2g-1)\text{ where }
\mathcal{L}^K=\bigcup_{(\gamma,\alpha)\in\mathcal{M}'_0} \mathcal{L}_{\gamma,\alpha}^K.$$
In fact, $(\alpha,\beta,\gamma,\kappa)$ can be served as coordinates of $\mathcal{L}^K$ (see \emph{p.236} in \cite{kri02}). We will call them the \emph{Krichever-Tyurin parameters}. Because of the dimension differences, the space $\mathcal{L}^{K}$ cannot be identified with a cotangent bundle $\TT^{\ast}\mathcal{M}'_0$ whose dimension is $2l^2g$. However, as in \cite{kri02} we may see that $\mathcal{L}^K/\mathbf{SL}(l,\mathbb{C})$ can be identified with 
$$\TT^\ast\widehat{\mathcal{M}'_0}=\TT^\ast\mathcal{M}'_0/\mathbf{SL}(l,\mathbb{C})\text{ where }\dim_\mathbb{C}\TT^\ast\widehat{\mathcal{M}'_0}=2(l^2(g-1)+1).$$
The Hitchin's abstract theory can be concretely realized by the Krichever-Lax matrices: Let $\ml(p;\gamma,\alpha)$ be an $(l\times l)$-Krichever-Lax-matrix on $\mathfrak{R}$ associated with Tyurin parameters $(\gamma,\alpha)$ and a canonical divisor $K$ where $\gamma=\gamma_1+\cdots+\gamma_{lg}$. Take a characteristic polynomial
$$R(\mu,p)=\det\Big(\mu\cdot I_{l\times l}-\ml(p;\gamma,\alpha)\Big)=0.$$
The zero locus $\{R(\mu,p)=0\}$ defines an algebraic curve. 
We denote the smoothly compactified model of this algebraic curve by $\widehat{\mathfrak{R}}$ and call it a \textit{spectral curve} associated with a Lax matrix $\ml(p;\gamma,\alpha)$. The coefficients $h_d(p;\ml)$ of
$$R(\mu,p)=\mu^l+\sum_{d=1}^{l}h_{d}(p;\ml)\mu^{l-d}$$
are a priori meromorphic functions on $\mathfrak{R}$ on the neighborhoods $U_j$ of $\gamma_j$ 
by definition. But it turns out that they are holomorphic on $U_j$ (see \emph{p.241} in \cite{kri02}). Now, we may have a map 
$$\mathbf{H}:\mathcal{L}^K\to\mathcal{H}^K\text{ by }\mathbf{H}(\ml)=\Big(h_1(p;\ml),\dots,h_l(p;\ml)\Big).$$
Since it is invariant under the conjugation action of $\mathbf{SL}(l,\mathbb{C})$, the map $\mathbf{H}$ can descend to the quotient space

\begin{equation}\label{hiteq1}
\mathbf{H}:\mathcal{L}^K/\mathbf{SL}(l,\mathbb{C})\to\mathcal{H}^K\text{ by }\mathbf{H}([\ml])
=\Big(h_1(p;\ml),\dots,h_l(p;\ml)\Big).
\end{equation}
This map is what Hitchin investigated in \cite{hit87}. By the parameters of the images and fibers of $\mathbf{H}$, we can foliate the space $\mathcal{L}^K/\mathbf{SL}(l,\mathbb{C})$. We summarize the contents of 
\emph{pp.241--243} in \cite{kri02} as follows:

\begin{theorem}\cite{kri02}\label{laxth1}
Let $[\ml]\in\mathcal{L}^K/\mathbf{SL}(l,\mathbb{C})$ be an $\mathbf{SL}(l,\mathbb{C})$-orbit of $\ml$ in $\mathcal{L}^K$. Then there is a one-to-one correspondence
$$[\ml]\longleftrightarrow\Big((h_1,\dots,h_l),[\widehat{D}]\Big)=\Big(\widehat{\mathfrak{R}},[\widehat{D}]\Big).$$
$[\widehat{D}]$ is an equivalence class of an effective divisor of degree $\widehat{g}+l-1$ on 
$\widehat{\mathfrak{R}}$ where $\widehat{g}=l^2(g-1)+1$ is the genus of $\widehat{\mathfrak{R}}$.
\end{theorem}

\section{Eigenvector mappings and the Euler sequence}
Let $\widehat{\mathfrak{R}}$ be a spectral curve associated with a Lax matrix $\ml(p;\gamma,\alpha)$:
$$\widehat{\mathfrak{R}}=\{\det\big(\mu\cdot I_{l\times l}-\ml(p;\gamma,\alpha)\big)=0\}\text{ where }p\in\mathfrak{R}.$$
Each point $(\mu,p):=\widehat{p}\in\widehat{\mathfrak{R}}$ is an eigenvalue of $\ml(p;\gamma,\alpha)$. From Theorem~\ref{laxth1}, it is not hard to see that for a Lax matrix $\ml(p;\gamma,\alpha)$, there exists a unique eigenspace complex line bundle $\LL$ of $\ml(p;\gamma,\alpha)$ on $\widehat{\mathfrak{R}}$ which is a sub-bundle of a trivial bundle $\mathbb{C}^l$ on $\widehat{\mathfrak{R}}$.

\begin{definition}\label{shde1}
We shall call \eqref{eqco1}
\begin{equation}\label{eqco1}
\overline{\boldsymbol{\psi}}_t(\gamma(t),\alpha(t)):\widehat{\mathfrak{R}}\to\mathbb{P}^{l-1}
\end{equation}
an \emph{eigenvector mapping} associated to Lax representation~\eqref{de5lax}.
\end{definition}

In other words, letting 
$\overline{\boldsymbol{\psi}}_t(\widehat{p};\gamma(t),\alpha(t))=\mathbb{C}\cdot\boldsymbol{\psi}_t(\widehat{p};\gamma(t),\alpha(t))$, we have
$$\boldsymbol{\psi}_t(\widehat{p};\gamma(t),\alpha(t))\ml_t(p;\gamma(t),\alpha(t))
=\mu(\widehat{p})\cdot\boldsymbol{\psi}_t(\widehat{p};\gamma(t),\alpha(t)).$$
A vector-valued meromorphic function $\boldsymbol{\psi}_t$ on $\widehat{\mathfrak{R}}$ defines a vector-valued (and multi-valued) meromorphic function $\pi_\ast\boldsymbol{\psi}_t$ on $\mathfrak{R}$ where $\pi:\widehat{\mathfrak{R}}\to\mathfrak{R}$ and the number of poles in each component of a vector $\boldsymbol{\psi}_t$ is $\widehat{g}+l-1$. Moreover, the multi-valued function $\pi_\ast\boldsymbol{\psi}_t$ has poles at $\gamma(t)=\gamma_1(t)+\cdots+\gamma_{lg}(t)$, and it is written as \eqref{kriform1} associated with a Tyurin parameter $(\gamma(t),\alpha(t))$ (see also Equation~\eqref{krieq2}). The eigenvalue $\mu(\widehat{p})$ can be regarded as a \emph{multi-valued} meromorphic function on $\mathfrak{R}$ with poles at the canonical divisor $K$ of $\mathfrak{R}$. Let 
$$\LL_{t}=\overline{\boldsymbol{\psi}}_{t}^{\ast}\big(\mathcal{O}_{\mathbb{P}^{l-1}}(1)\big)\in
\pic^{\widehat{g}+l-1}(\widehat{\mathfrak{R}}).$$
Note that the degree of $\LL_t$ is $\widehat{g}+l-1$ and $\LL_t$ is a line bundle associated with an equivalence class $[\widehat{D}_t]$ of divisors by Theorem~\ref{laxth1}. Let $\mathcal{L}^{K}_{\widehat{\mathfrak{R}}}/\mathbf{SL}(l,\mathbb{C})\subset \mathcal{L}^K/\mathbf{SL}(l,\mathbb{C})$ be the pre-images of Hitchin map~\eqref{hiteq1} associated to a spectral curve $\widehat{\mathfrak{R}}$. The eigenvector mapping $\overline{\boldsymbol{\psi}}_t$ induces 

\begin{equation}\label{eqco2}
\varphi_{\widehat{\mathfrak{R}}}:\mathcal{L}^{K}_{\widehat{\mathfrak{R}}}/\mathbf{SL}(l,\mathbb{C})\to\pic^{\widehat{g}+l-1}(\widehat{\mathfrak{R}})\text{ by }\varphi_{\widehat{\mathfrak{R}}}([\ml_t])=\overline{\boldsymbol{\psi}}_{t}^{\ast}\big(\mathcal{O}_{\mathbb{P}^{l-1}}(1)\big).
\end{equation}
We will also call this map $\varphi_{\widehat{\mathfrak{R}}}$ an \emph{eigenvector mapping} associated to a spectral curve $\widehat{\mathfrak{R}}$. Since the tangent space of $\pic^{\widehat{g}+l-1}(\widehat{\mathfrak{R}})$ is isomorphic to $\HH^1(\widehat{\mathfrak{R}},\mathcal{O}_{\widehat{\mathfrak{R}}})$, we have
$$\frac{d}{dt}\LL_t|_{t=0}\in \HH^1(\widehat{\mathfrak{R}},\mathcal{O}_{\widehat{\mathfrak{R}}}).$$

Consider the \emph{Euler sequence} over $\mathbb{P}^{l-1}$
\begin{equation}\label{euler}
\xymatrix{
0\ar[r]&\mathcal{O}_{\mathbb{P}^{l-1}}\ar[r]&\mathbb{C}^l\otimes \mathcal{O}_{\mathbb{P}^{l-1}}(1)\ar[r]&\varTheta_{\mathbb{P}^{l-1}} \ar[r]&0.}
\end{equation}
Here $\varTheta_{\mathbb{P}^{l-1}}$ is the sheaf of a \emph{holomorphic} tangent bundle $\TT\mathbb{P}^{l-1}$. Pulling back Euler sequence~\eqref{euler} to $\widehat{\mathfrak{R}}$ by $\overline{\boldsymbol{\psi}}_t$ induces the following short exact sequence on  $\widehat{\mathfrak{R}}$:

\begin{equation}\label{exac1}
\xymatrix{
0\ar[r]&\mathcal{O}_{\widehat{\mathfrak{R}}}\ar[r]&\mathbb{C}^l\otimes \LL_{t}\ar[r]&\overline{\boldsymbol{\psi}}_{t}^{\ast}\varTheta_{\mathbb{P}^{l-1}} \ar[r]&0.}
\end{equation}
From short exact sequence~\eqref{exac1}, we have a long exact sequence:

\begin{equation}\label{lexac1}
\xymatrix{
\cdots\ar[r]&\HH^0(\widehat{\mathfrak{R}},\mathbb{C}^l\otimes \LL_{t})\ar[r]&
\HH^0(\widehat{\mathfrak{R}},\overline{\boldsymbol{\psi}}_{t}^{\ast}\varTheta_{\mathbb{P}^{l-1}} )\ar[r]^(.55){\delta}&
\HH^1(\widehat{\mathfrak{R}},\mathcal{O}_{\widehat{\mathfrak{R}}})\ar[r]&\cdots.}
\end{equation}
Since $\overline{\boldsymbol{\psi}}_t=\mathbb{C}\cdot\boldsymbol{\psi}_t$, any global section $\mathbf{s}_t$ of $\mathbb{C}^l\otimes\LL_{t}\cong\bigoplus^l\LL_{t}$ can be given by $\{\rho_{t,i}^{-1}(\widehat{p})\cdot\boldsymbol{\psi}_{t,i}\}$ where $\boldsymbol{\psi}_{t,i}$ is the restriction of 
$\boldsymbol{\psi}_t$ to an open set $U_i$ and 
$$\rho_{t,i}(\widehat{p})^{-1}\cdot \rho_{t,j}(\widehat{p})=g_{t,ij}(\widehat{p})\text{ on }U_i\cap U_j.$$
Here $\{\rho_{t,i}(\widehat{p})\}$ is the set of local non-vanishing holomorphic functions.

\begin{lemma}\label{cohomle1}
A time-derivative $\frac{d}{dt}\mathbf{s}_t$ can be regarded as an element of $\HH^0(\widehat{\mathfrak{R}},\overline{\boldsymbol{\psi}}_{t}^{\ast}\varTheta_{\mathbb{P}^{l-1}} )$.
\end{lemma}

\begin{proof}
Let $\mathbf{s}_t:=\{\rho_{t,i}^{-1}(\widehat{p})\cdot\boldsymbol{\psi}_{t,i}\}$ be a global section of  $\mathbb{C}^l\otimes\LL_{t}$. Accordingly, 
$$\frac{d}{dt}\boldsymbol{\psi}_{t,i}=\frac{d}{dt}(\rho_{t,i}\cdot\mathbf{s}_{t,i})=(\frac{d}{dt}\rho_{t,i})\cdot\mathbf{s}_{t,i}+\rho_{t,i}\cdot(\frac{d}{dt}\mathbf{s}_{t,i}).$$
Here $\mathbf{s}_{t,i}$ is the restriction of $\mathbf{s}_{t}$ to $U_i$. Since
$$\xymatrix@1@R=1pt{
0\ar[r]&\mathcal{O}_{\widehat{\mathfrak{R}}}\ar[r]^(.4){\mathbf{s}_t}&\mathbb{C}^l\otimes \LL_{t}
\ar[r]& \overline{\boldsymbol{\psi}}_{t}^{\ast}\varTheta_{\mathbb{P}^{l-1}} \ar[r]&0,\\
&\{\xi_i\}\ar[r]& \{\xi_i\cdot\mathbf{s}_{t,i}\}&&&}$$
$\{\rho_{t,i}^{-1}\cdot\frac{d}{dt}\boldsymbol{\psi}_{t,i}\}=[\frac{d}{dt}\mathbf{s}_t]$ defines an element of $\HH^0(\widehat{\mathfrak{R}},\mathbb{C}^l\otimes \LL_{t}/\mathcal{O}_{\widehat{\mathfrak{R}}})$. From 
$$\HH^0(\widehat{\mathfrak{R}},\mathbb{C}^l\otimes \LL_{t}/\mathcal{O}_{\widehat{\mathfrak{R}}})\cong \HH^0(\widehat{\mathfrak{R}},\overline{\boldsymbol{\psi}}_{t}^{\ast}\varTheta_{\mathbb{P}^{l-1}} ),$$
we can see that  $\{\rho_{t,i}^{-1}\cdot\frac{d}{dt}\boldsymbol{\psi}_{t,i}\}$ defines an element of $\HH^0(\widehat{\mathfrak{R}},\overline{\boldsymbol{\psi}}_{t}^{\ast}\varTheta_{\mathbb{P}^{l-1}} )$.
\end{proof}

The mapping $\varphi_{\widehat{\mathfrak{R}}}:\mathcal{L}^{K}_{\widehat{\mathfrak{R}}}/\mathbf{SL}(l,\mathbb{C})\to\pic^{\widehat{g}+l-1}(\widehat{\mathfrak{R}})$ induces a mapping between tangent spaces
$$\TT\varphi_{\widehat{\mathfrak{R}}}:\TT_{[\ml]}\mathcal{L}^{K}_{\widehat{\mathfrak{R}}}/\mathbf{SL}(l,\mathbb{C})\to\HH^1(\widehat{\mathfrak{R}},\mathcal{O}_{\widehat{\mathfrak{R}}})\text{ where }[\ml]\in\mathcal{L}^{K}_{\widehat{\mathfrak{R}}}/\mathbf{SL}(l,\mathbb{C}).$$
In other words,
$$\TT\varphi_{\widehat{\mathfrak{R}}}
\big(\frac{d}{dt}[\ml_t(p;\gamma(t),\alpha(t))]|_{t=0}\big)\in\HH^1(\widehat{\mathfrak{R}},\mathcal{O}_{\widehat{\mathfrak{R}}}).$$
We can observe the following result:

\begin{theorem}\label{cohomth1}
Let $\mathbf{s}_t:=\{\rho_{t,i}^{-1}(\widehat{p})\cdot\boldsymbol{\psi}_{t,i}\}$ be a global section of  $\mathbb{C}^l\otimes\LL_{t}$ and $[\frac{d}{dt}\mathbf{s}_t]$ in Lemma~\ref{cohomle1} be regarded as 
an element of  $\HH^0(\widehat{\mathfrak{R}},\overline{\boldsymbol{\psi}}_{t}^{\ast}\varTheta_{\mathbb{P}^{l-1}} )$. Then we may have
$$\TT\varphi_{\widehat{\mathfrak{R}}}\big(\frac{d}{dt}[\ml_t]|_{t=0}\big)
=\delta([\frac{d}{dt}\mathbf{s}_t|_{t=0}]).$$
Moreover, it is independent of a section $\mathbf{s}_t$ we chose.
\end{theorem}

\begin{proof}
The infinitesimal change $\TT\varphi_{\widehat{\mathfrak{R}}}
\big(\frac{d}{dt}[\ml_t]|_{t=0}\big)$ of line bundles 
is characterized by $\frac{d}{dt}\log g_{ij}(t)$ where $\{g_{ij}(t)\}$ is the set of transition functions of a line bundle $\LL_t$ over $\widehat{\mathfrak{R}}$ associated with an open cover $\{U_i\}$. The connecting homomorphism $\delta$ of long exact sequence~\eqref{lexac1} is given by 
$$\delta([\frac{d}{dt}\mathbf{s}_t])=\{\rho_{t,j}^{-1}\cdot(\frac{d}{dt}\rho_{t,j})-\rho_{t,i}^{-1}\cdot(\frac{d}{dt}\rho_{t,i})\}\in\HH^1(\widehat{\mathfrak{R}},\mathcal{O}_{\widehat{\mathfrak{R}}}).$$
Since $\rho_{t,i}^{-1}\cdot \rho_{t,j}=g_{ij}(t)$ on $U_i\cap U_j$, we have
$$\rho_{t,j}^{-1}\cdot(\frac{d}{dt}\rho_{t,j})-\rho_{t,i}^{-1}\cdot(\frac{d}{dt}\rho_{t,i})=\frac{d}{dt}\log g_{ij}(t).$$
So, $\TT\varphi_{\widehat{\mathfrak{R}}}\big(\frac{d}{dt}\ml_t|_{t=0}\big)
=\delta([\frac{d}{dt}\mathbf{s}_t|_{t=0}])$. Also notice that the explicit expression 
$$\delta([\frac{d}{dt}\mathbf{s}_t])=\{\rho_{t,j}^{-1}\cdot(\frac{d}{dt}\rho_{t,j})-\rho_{t,i}^{-1}\cdot(\frac{d}{dt}\rho_{t,i})\}\in\HH^1(\widehat{\mathfrak{R}},\mathcal{O}_{\widehat{\mathfrak{R}}})$$
implies that it only depends on $\boldsymbol{\psi}_t$ and is independent of choosing $\mathbf{s}_t$.
\end{proof}

\section{A cohomological interpretation of straight line flows in the space of Krichever-Lax matrices}\label{resec3}

The dynamics of Lax representation~\eqref{de5lax} on $\mathcal{L}^K$ is invariant under the addition of a polynomial $\mathsf{P}(\ml_t)$ of $\ml_t$ or an element $\mq_t$ commuting with $\ml_t$, since 

\begin{equation}\label{co49}
\frac{d}{dt}\ml_t=[\mm_t,\ml_t]=[\mm_t+\mathsf{P}(\ml_t),\ml_t]=[\mm_t+\mq_t,\ml_t].
\end{equation}
Thus, the dependence of flows on $\mm$ might be indicated by an equivalence object associated with $\mm$. We will characterize it in terms of a cohomological class associated with $\mm$. In fact, what we are interested in is a flow in the quotient space $\mathcal{L}^K/\mathbf{SL}(l,\mathbb{C})$. Note that for $W\in\mathbf{SL}(l,\mathbb{C})$, we have

\begin{equation}\label{ineq1}
\begin{aligned}
\frac{d}{dt}(W^{-1}\ml W)&=[\mm,W^{-1}\ml W]=\mm(W^{-1}\ml W)-(W^{-1}\ml W)\mm\\
&=W^{-1}(W\mm W^{-1}\ml -\ml W\mm W^{-1})W\\
&=W^{-1}[W\mm W^{-1},\ml]W.
\end{aligned}
\end{equation}
Thus, if $\frac{d}{dt}\ml=[\mm,\ml]$, 
then $\frac{d}{dt}\ml=[W\mm W^{-1},\ml]$. So, the characteristic class of $\mm$ should be invariant under the change of gauges. We will show the gauge-invariance of the associated cohomology class of $\mm_t$ in Lemma~\ref{lem4}.

First, we describe the condition on isospectral deformations, that is, the condition that the flow of a Lax representation stays in a leaf $\mathcal{L}^{K}_{\widehat{\mathfrak{R}}}$ in the foliation of the Hitchin map.

\begin{lemma}\label{krle3}
If the flow of a vector field $[\mm_t,\ml_t]$ is tangent to 
$\mathcal{L}^{K}$, then $[\mm_t,\ml_t]$ has poles only at the canonical divisor $K$ of $\mathfrak{R}$ other than $\gamma(t)=\gamma_1(t)+\cdots+\gamma_{lg}(t)$.
\end{lemma}

\begin{proof}
Suppose that a vector field $[\mm_t,\ml_t]$ on the space of matrix-valued meromorphic functions on $\mathfrak{R}$ is tangent to $\mathcal{L}^{K}$. Then the flow $\ml_t$ stays in 
$\mathcal{L}^{K}$. So, we can write $\frac{d}{dt}\ml_t=[\mm_t,\ml_t]$. From Definition~\ref{laxde1}, it is easy to see that $\frac{d}{dt}\ml_t$ has a double pole possibly at $\gamma_{j}(t)$ for $j=1,\dots,lg$ and a simple pole at $p_i$ where $K=\sum_{i=1}^{2g-2}p_i$. Thus, we have the desired result.
\end{proof}

Suppose that $[\mm_t,\ml_t]$ is tangent to $\mathcal{L}^{K}_{\widehat{\mathfrak{R}}}$. From eigenvector mapping~\eqref{eqco2}, we have $\boldsymbol{\psi}_t\ml_t=\mu\cdot\boldsymbol{\psi}_t$. After differentiating $\boldsymbol{\psi}_t\ml_t=\mu\cdot\boldsymbol{\psi}_t$ with respect to $t$, we have $(\frac{d}{dt}\boldsymbol{\psi}_t)\ml_t+\boldsymbol{\psi}_t(\frac{d}{dt}\ml_t)=\mu\cdot(\frac{d}{dt}\boldsymbol{\psi}_t)$. Note that $\mu(\widehat{p})$ does not depend on $t$, i.e., it is isospectral. From $\boldsymbol{\psi}_t\ml_t=\mu\cdot\boldsymbol{\psi}_t$ and 
$(\frac{d}{dt}\boldsymbol{\psi}_t)\ml_t+\boldsymbol{\psi}_t[\mm_t,\ml_t]=\mu\cdot(\frac{d}{dt}\boldsymbol{\psi}_t)$, we have
$$\big(\boldsymbol{\psi}_t\mm_t+(\frac{d}{dt}\boldsymbol{\psi}_t)\big)\ml_t=
\mu\cdot\big(\boldsymbol{\psi}_t\mm_t+(\frac{d}{dt}\boldsymbol{\psi}_t)\big).$$
Since the eigenspace of $\ml_t(p)$ associated with the eigenvalue $\mu$ is $1$-dimensional generically, 
we find a meromorphic function $\lambda_t(\widehat{p};\gamma(t),\alpha(t))$ such that

\begin{equation}\label{coseq1}
\boldsymbol{\psi}_t\mm_t+(\frac{d}{dt}\boldsymbol{\psi}_t)=\lambda_t\boldsymbol{\psi}_t.
\end{equation}
This meromorphic function $\lambda_t$ certainly depends on $\mm_t$ and $\boldsymbol{\psi}_t$. However, the Laurent tails of $\lambda_t$ at poles only depend on $\mm_t$: For another $\boldsymbol{\psi}'_t=\varrho_t\cdot\boldsymbol{\psi}_t$ associated with a line bundle $\LL_{t}$ where $\varrho_t(\widehat{p})$ is a local non-vanishing holomorphic function on $\widehat{\mathfrak{R}}$, $\lambda_t$ is transformed to $\lambda_t+\varrho_{t}^{-1}\frac{d}{dt}\varrho_t$. Thus the Laurent tails are well-defined quantities associated to $\mm_t$ only. 
Hence, the meromorphic functions $\lambda_t$ in Equation~\eqref{coseq1} can be regarded as a global section of a skyscraper sheaf $\mathbb{C}_{\pi^{-1}(nK)}$ for some positive integer $n$ where $\pi:\widehat{\mathfrak{R}}\to\mathfrak{R}$ and $K$ is a canonical divisor of $\mathfrak{R}$. In this notation we make a definition:

\begin{definition}\label{rede2}
A \emph{residue section} $\rho(\mm_t)\in \HH^0(\widehat{\mathfrak{R}},\mathbb{C}_{\pi^{-1}(nK)})$ associated to 
$\mm_t$ is defined to be the \emph{Laurent tail} $\{\lambda_{t,i}\}$ of $\lambda_t$ in Equation~\eqref{coseq1} at $\pi^{-1}(nK)$ where $K=\sum_{i=1}^{2g-2}p_i$.
\end{definition}

Generically, the poles of $\boldsymbol{\psi}_t$ are simple and $\deg(\boldsymbol{\psi}_t)_\infty=\widehat{g}+l-1$. So in the neighborhood of a pole $\widehat{\gamma}_{j}(t)$ of $\boldsymbol{\psi}_t$ we may write $\boldsymbol{\psi}_t$ as
$$\boldsymbol{\psi}_t(\widehat{z})=\frac{\mathbf{c}_j(t)}{\widehat{z}-\widehat{z}(\widehat{\gamma}_j(t))}+O(1)\text{ where }\mathbf{c}_j(t)\in\mathbb{C}^l.$$
Consequently,

\begin{equation}\label{ncoeq1}
\frac{d}{dt}\boldsymbol{\psi}_t=\frac{\mathbf{c}_j(t)\cdot\frac{d}{dt}\widehat{z}(\widehat{\gamma}_j(t))}{\big(\widehat{z}-\widehat{z}(\gamma_j(t))\big)^2}+\frac{\frac{d}{dt}\mathbf{c}_j(t)}{\widehat{z}-\widehat{z}(\gamma_j(t))}+O(1).
\end{equation}
The next theorem indicates how the behavior of the poles of the global meromorphic function $\lambda_t$ on $\widehat{\mathfrak{R}}$ governs the dynamics of Lax representation.

\begin{theorem}\label{thco1}
Suppose that $\frac{d}{dt}\ml_t=[\mm_t,\ml_t]$ and $(s_0)=\pi^{-1}(nK)$. There is 
$\lambda_t\cdot s_0\in \HH^0(\widehat{\mathfrak{R}},\pi^\ast\K_{\mathfrak{R}}^n)$ for some positive integer $n$ such that $\boldsymbol{\psi}_t\mm_t-\lambda_t\boldsymbol{\psi}_t$ defines a global section of $\mathbb{C}^{l}\otimes \LL_{t}$ if and only if the flows are constant, i.e., $\frac{d}{dt}\ml_t=0$.
\end{theorem}

\begin{proof}
Suppose that  there is 
$\lambda_t\cdot s_0\in \HH^0(\widehat{\mathfrak{R}},\pi^\ast\K_{\mathfrak{R}}^n)$ for some positive integer $n$ such that $\boldsymbol{\psi}_t\mm_t-\lambda_t\boldsymbol{\psi}_t$ defines a global section of $\mathbb{C}^{l}\otimes \LL_{t}$. Accordingly, there is a global meromorphic function $\xi_t$ on $\widehat{\mathfrak{R}}$ such that 
$$\boldsymbol{\psi}_t\mm_t-\lambda_t\boldsymbol{\psi}_t=\xi_t\cdot\boldsymbol{\psi}_t.$$
Of course, the only possible poles of $\xi_t$ are at $\pi^{-1}(nK)$, since $\lambda_t\cdot s_0\in \HH^0(\widehat{\mathfrak{R}},\pi^\ast\K_{\mathfrak{R}}^n)$.  This implies that $\mm_t$ preserves the eigenspaces of $\ml_t$. Thus, $\mm_t$ and $\ml_t$ commute. From $\frac{d}{dt}\ml_t=[\mm_t,\ml_t]$, we conclude that $\frac{d}{dt}\ml_t=0$.

Suppose that $\frac{d}{dt}\ml_t=0$. Since $[\mm_t,\ml_t]=0$, the $\mm_t$ preserves the eigenspaces of $\ml_t$. What this amounts is that there is a global meromorphic function $\varsigma_t(\widehat{p})$ on $\widehat{\mathfrak{R}}$ 
such that $\boldsymbol{\psi}_t\mm_t=\varsigma_t\boldsymbol{\psi}_t$. Notice that $\varsigma_t$ only has poles possibly at $\pi^{-1}(nK)$ where $n$ is a positive integer and $K$ is the canonical divisor of $\mathfrak{R}$, since $\mm$ preserves the eigenspace of $\ml$ and $\boldsymbol{\psi}_t\ml_t=\mu\cdot\boldsymbol{\psi}_t$ where $\mu(\widehat{p})$ takes poles only at $\pi^{-1}(K)$. Thus,
$$\boldsymbol{\psi}_t\mm_t-\lambda_t\boldsymbol{\psi}_t=(\varsigma_t-\lambda_t)\boldsymbol{\psi}_t$$
defines a global section of $\mathbb{C}^{l}\otimes \LL_{t}$. 
Moreover, since $\frac{d}{dt}\ml_t=0$ implies that $\frac{d}{dt}\widehat{z}(\widehat{\gamma}_j(t))=0$ for $j=1,\dots,\widehat{g}+l-1$, we see that $\frac{d}{dt}\boldsymbol{\psi}_t$ has only first order poles at $\gamma_j$ from Equations~\eqref{ncoeq1}. Since $(\varsigma_t-\lambda_t)\boldsymbol{\psi}_t=-\frac{d}{dt}\boldsymbol{\psi}_t$, we conclude that $\lambda_t\in \HH^0(\widehat{\mathfrak{R}},\pi^\ast\K_{\mathfrak{R}}^{n})$.
\end{proof}

Theorem~\ref{thco1} exhibits how the dynamics on $\mathcal{L}^{K}_{\widehat{\mathfrak{R}}}$ of Lax representation~\eqref{de5lax} is related with $\mm_t$ in terms of the residue section $\rho(\mm_t)=\{\lambda_{t,i}\}\in\HH^0(\widehat{\mathfrak{R}},\mathbb{C}_{\pi^{-1}(nK)})$. When the flow is constant, then $\rho(\mm_t)=\{\lambda_{t,i}\}$ defines a global section in $\HH^0(\widehat{\mathfrak{R}},\pi^\ast\K_{\mathfrak{R}}^n)$. In other words, if $\ml_t$ and $\mm_t$ commute, then $\mm_t$ defines an endomorphism of $\E_{\gamma(t),\alpha(t)}$. In Corollary~\ref{cohomco3}, we will give a necessary and sufficient condition for $\{\lambda_{t,i}\}$ when the flow of $\ml_t$ is linear, which is the second simplest case next to the constant flows.

It is not hard to see that a residue section $\rho(\mm_t)$ of $\mm_t$ is gauge-invariant:

\begin{lemma}\label{lem4}
For $W\in \mathbf{SL}(l,\mathbb{C})$, we have $\rho(\mm_t)=\rho(W^{-1}\cdot\mm_t\cdot W)$.
\end{lemma}

\begin{proof}
Let $\boldsymbol{\psi}_t\ml_t=\mu\cdot\boldsymbol{\psi}_t$ and 
$\boldsymbol{\psi}_t\mm_t+\frac{d}{dt}\boldsymbol{\psi}_t=\lambda_t\cdot\boldsymbol{\psi}_t$. Since
$$(W\boldsymbol{\psi}_tW^{-1})\ml_t=W\cdot\boldsymbol{\psi}_t(W^{-1}\ml_tW)\cdot W^{-1}=\mu\cdot (W\boldsymbol{\psi}_tW^{-1}),$$
we have $(W\boldsymbol{\psi}_tW^{-1})\mm_t+\frac{d}{dt}(W\boldsymbol{\psi}_tW^{-1})=\lambda_t\cdot W\boldsymbol{\psi}_tW^{-1}$. Accordingly,
$$\begin{aligned}
\boldsymbol{\psi}_t(W^{-1}\mm_t W)+\frac{d}{dt}\boldsymbol{\psi}_t&=W^{-1}\cdot\big((W\boldsymbol{\psi}_tW^{-1})\mm_t+\frac{d}{dt}(W\boldsymbol{\psi}_tW^{-1})\big)W\\
&=W^{-1}\cdot(\lambda_t\cdot W\boldsymbol{\psi}_tW^{-1})\cdot W=\lambda_t\boldsymbol{\psi}_t.
\end{aligned}$$
\end{proof}

For a positive integer $n$, consider the short exact sequence of elementary transformation

\begin{equation}\label{exsk1}
\xymatrix{
0\ar[r]&\mathcal{O}_{\widehat{\mathfrak{R}}}\ar[r]&\mathcal{O}_{\widehat{\mathfrak{R}}}\otimes\pi^\ast \K^n\ar[r]^{\jmath}&\mathbb{C}_{\pi^{-1}(nK)}\ar[r]&0.}
\end{equation}
This induces a long exact sequence

\begin{equation}\label{lexsk1}
\xymatrix@R=1pt{
0\ar[r]&\HH^0(\widehat{\mathfrak{R}},\mathcal{O}_{\widehat{\mathfrak{R}}})\ar[r]&\HH^0(\widehat{\mathfrak{R}},\pi^\ast\K^n)\ar[r]^(.45){\jmath}&\HH^0(\widehat{\mathfrak{R}},\mathbb{C}_{\pi^{-1}(nK)})&\\
\ar[r]^(.3){\partial}&\HH^1(\widehat{\mathfrak{R}},\mathcal{O}_{\widehat{\mathfrak{R}}})\ar[r]&\HH^1(\widehat{\mathfrak{R}},\pi^\ast \K^n)\ar[r]&\HH^1(\widehat{\mathfrak{R}},\mathbb{C}_{\pi^{-1}(nK)}).&}
\end{equation}
The time dependence of the residue section $\rho(\mm_t)=\{\lambda_{t,i}\}$ associated to $\mm_t$ can be characterized by the following theorem:

\begin{theorem}\label{cothe1}
Let $\TT\varphi_{\widehat{\mathfrak{R}}}\big(\frac{d}{dt}[\ml_t]\big)=\frac{d}{dt}\LL_t$. If $[\mm_t,\ml_t]$ is tangent to $\mathcal{L}^{K}_{\widehat{\mathfrak{R}}}$, then   

\begin{equation}\label{cotheq1}
\frac{d}{dt}\LL_t=\partial\rho(\mm_t) \in\HH^1(\widehat{\mathfrak{R}},\mathcal{O}_{\widehat{\mathfrak{R}}})\cong\HH^0(\widehat{\mathfrak{R}},\K_{\widehat{\mathfrak{R}}}).
\end{equation}
\end{theorem}

\begin{proof}
We let $\varpi_1=\{\rho_{t,i}^{-1}\cdot\frac{d}{dt}\boldsymbol{\psi}_t\}\in \HH^0(\widehat{\mathfrak{R}},\overline{\boldsymbol{\psi}}_{t}^{\ast}\varTheta_{\mathbb{P}^{l-1}})$ in the notation of the proof of Lemma~\ref{cohomle1}. Similarly, we may let
$$\varpi_2=\{\rho_{t,i}^{-1}\cdot\lambda_{t,i}\cdot\boldsymbol{\psi}_t\}=\{\rho_{t,i}^{-1}\cdot(\boldsymbol{\psi}_t\mm_t+\frac{d}{dt}\boldsymbol{\psi}_t)\}\in \HH^0(\widehat{\mathfrak{R}},\mathbb{C}^l\otimes \LL_{t}\otimes\pi^\ast\K^n).$$
Since $\HH^0(\widehat{\mathfrak{R}},\mathbb{C}^l\otimes \LL_{t}\otimes\pi^\ast\K^n/
\mathcal{O}_{\widehat{\mathfrak{R}}})\cong \HH^0(\widehat{\mathfrak{R}},\overline{\boldsymbol{\psi}}_{t}^{\ast}\varTheta_{\mathbb{P}^{l-1}}\otimes\pi^\ast\K^n )$, $\varpi_2$ may induce an element in $\HH^0(\widehat{\mathfrak{R}},\overline{\boldsymbol{\psi}}_{t}^{\ast}\varTheta_{\mathbb{P}^{l-1}}\otimes\pi^\ast\K^n )$. Let us denote this element by $\tau(\varpi_2)$. Now we let $\varpi_3=\{\lambda_{t,i}\}\in \HH^0(\widehat{\mathfrak{R}},\mathbb{C}_{\pi^{-1}(nK)})$. Since $\rho_{t,i}$ is a non-vanishing local holomorphic function,  
from 
$$\xymatrix@R=1pt{
0\ar[r]&\mathcal{O}_{\widehat{\mathfrak{R}}}\ar[r]&\mathcal{O}_{\widehat{\mathfrak{R}}}\otimes\pi^\ast \K^n\ar[r]^{\jmath}&\mathbb{C}_{\pi^{-1}(nK)}\ar[r]&0,\\
&&
\rho_{t,i}^{-1}\cdot\lambda_{t,i}\ar[r]^(.6){\jmath}&\lambda_{t,i}&}$$
it is clear that $\partial\{\lambda_{t,i}\}=\{\frac{d}{dt}\log g_{t,ij}\}$. From Theorem~\ref{cohomth1}, we also have
$$\delta\{\rho_{t,i}^{-1}\cdot\frac{d}{dt}\boldsymbol{\psi}_t\}=\{\rho_{t,j}^{-1}\cdot\frac{d}{dt}\rho_{t,j}-\rho_{t,i}^{-1}\cdot\frac{d}{dt}\rho_{t,i}\}=\{\frac{d}{dt}\log g_{t,ij}\}.$$
Hence, $\frac{d}{dt}\LL_t=\partial\rho(\mm_t)$. Cohomologically, this is just chasing the following diagram:

$$\small\xymatrix{
\HH^0(\widehat{\mathfrak{R}},\mathcal{O}_{\widehat{\mathfrak{R}}})\ar[r]\ar[d]&\HH^0(\widehat{\mathfrak{R}},\pi^\ast\K^n)\ar[d]\ar[r]^(.45){\jmath}&\HH^0(\widehat{\mathfrak{R}},\mathbb{C}_{\pi^{-1}(nK)})\ar[d]^{\sigma}\ar[r]^(0.7){\partial}&\\
\HH^0(\widehat{\mathfrak{R}},\mathbb{C}^l\otimes \LL_{t})\ar[d]\ar[r]&\HH^0(\widehat{\mathfrak{R}},\mathbb{C}^l\otimes \LL_{t}\otimes\pi^\ast \K^n)\ar[d]^{\tau}\ar[r]^(.45){\jmath}&\HH^0(\widehat{\mathfrak{R}},\mathbb{C}^l\otimes \LL_{t}\otimes\mathbb{C}_{\pi^{-1}(nK)})\ar[d]^{\tau}&\\
\HH^0(\widehat{\mathfrak{R}},\overline{\boldsymbol{\psi}}_{t}^{\ast}\varTheta_{\mathbb{P}^{l-1}} )\ar[d]^{\delta}\ar[r]^(.4){\imath}&
\HH^0(\widehat{\mathfrak{R}},\overline{\boldsymbol{\psi}}_{t}^{\ast}\varTheta_{\mathbb{P}^{l-1}}\otimes\pi^\ast\K^n)\ar[r]^(.45){\jmath}&
\HH^0(\widehat{\mathfrak{R}},\overline{\boldsymbol{\psi}}_{t}^{\ast}\varTheta_{\mathbb{P}^{l-1}} \otimes\mathbb{C}_{\pi^{-1}(nK)})&\\
\HH^1(\widehat{\mathfrak{R}},\mathcal{O}_{\widehat{\mathfrak{R}}})&&&}$$
Since $\tau\circ\jmath=\jmath\circ\tau$ and $\imath(\varpi_1)=\tau(\varpi_2)$, we have 
$$\tau\circ\jmath(\varpi_2)=\jmath\circ\tau(\varpi_2)=\jmath\circ\imath(\varpi_1)=0.$$
Hence, there is $\varpi_3\in \HH^0(\widehat{\mathfrak{R}},\mathbb{C}_{\pi^{-1}(nK)})$ such that $\sigma(\varpi_3)=\jmath(\varpi_2)$. From the chasing the diagram, we see $\partial(\varpi_3)=\delta(\varpi_1)$.
\end{proof}

Note that Theorem~\ref{cothe1} confirms Theorem~\ref{thco1} again and this cohomological proof of Theorem~\ref{cothe1} again shows that the gauge-invariance of the residue section, which was verified in Lemma~\ref{lem4}. Moreover, we can deduce from Theorem~\ref{cothe1} that the flow on the quotient of $\mathcal{L}^{K}_{\widehat{\mathfrak{R}}}$ is described by the Laurent tails of $\lambda_t$ at $\pi^{-1}(nK)$ of $\widehat{\mathfrak{R}}$. A corollary we can have from Theorem~\ref{cothe1} is as follows:

\begin{corollary}\label{cohomco3}
$\LL_t$ is linear on $\pic^{\widehat{g}+l-1}(\widehat{\mathfrak{R}})$ if and only if 
$$\frac{d}{dt}\rho(\mm_t)\equiv0\textnormal{ modulo }\spa\{\jmath\big(\HH^0(\widehat{\mathfrak{R}},\pi^\ast \K^n)\big) ,\rho(\mm_t)\}.$$
\end{corollary}

\begin{proof}
Clearly, we can observe that the flow $\LL_t$ is straight if $\frac{d^2}{dt^2}\LL_t=0$ or $\frac{d^2}{dt^2}\LL_t=c\cdot\frac{d}{dt}\LL_t$ where $c\ne0$. By Theorem~\ref{cothe1}, $\frac{d^2}{dt^2}\LL_t=0$ if and only if $\frac{d}{dt}\rho(\mm_t)\equiv0$ modulo $\jmath\big(\HH^0(\widehat{\mathfrak{R}},\pi^\ast \K^n)\big)$. And $\frac{d^2}{dt^2}\LL_t=c\cdot\frac{d}{dt}\LL_t$ if and only if $\frac{d}{dt}\rho(\mm_t)\equiv0$ modulo $\rho(\mm_t)$. This proves the claim.
\end{proof}

\section{A characterization of flows in terms of $\mm$}\label{secch4}

As in the case of meromorphic functions on a compact Riemann surface, a matrix-valued meromorphic function on a compact Riemann surface is determined by the behavior of its poles. Consequently, the characterization of poles of $\mm_t$ determines the dynamics of Lax representation~\eqref{de5lax}. From Lemma~\ref{krle3}, we may see that at the poles of $\mm_t$ other than $lg$ points $\gamma_j$, the poles of $[\mm_t,\ml_t]$ are no greater than the poles of $\ml_t$. This is one restriction for defining tangent flows and it turns out to be the only one.

The existence of a meromorphic (matrix-valued) function on a compact Riemann surface is manifested by the Riemann-Roch theorem. Accordingly, we may not have $\mm$ for generally prescribed poles $D$. What this means is that we need special ansatz to have the existence of $\mm$. 
In \emph{p.233} of \cite{kri02}, Krichever defines special ansatz which guarantees the existence of $\mm$. That is, $\mm$ exits if $\mm$ has a special form of Equation~\eqref{lemmaeq1} at $lg$ points 
\begin{equation}\label{lemmaeq1}
\mm=\frac{\mm_{j,-1}}{z-z(\gamma_j)}+\mm_{j,0}+O(z-z(\gamma_j))\text{ for }j=1,\dots,lg,
\end{equation}
where the $(l\times l)$-matrix $\mm_{j,-1}$ 
is given by $\mathbf{v}^{T}_{j}\cdot\boldsymbol{\alpha}_{j}$ for $\mathbf{v}_j\in\mathbb{C}^l$.

The description of straight line flows in terms of $\mm_t$ will be given as follows: Let $K=\sum_{i=1}^{2g-2}p_i$ be a canonical divisor of $\mathfrak{R}$ where all $p_i$ are distinct. Consider $\mm_t$ satisfying Equation~\eqref{lemmaeq1} around $\gamma_j$ for 
$j=1,\dots,lg$ and locally given by 

\begin{equation}\label{chm}
\mm_t(w_i)=w^{-{m_i}}_{i}\ml^{n_i}_{t}\text{ around }p_i. 
\end{equation}
Here $w_i$ is a local coordinate around $p_i$. From Lemma~\ref{krle3}, we see that $[\mm_t,\ml_t]$ is tangent to $\mathcal{L}_{\widehat{\mathfrak{R}}}^{K}$. 
By Equation~\eqref{krieq22} (see \cite{kim07} for details)
\begin{equation}\label{krieq22}
\pi_\ast\boldsymbol{\psi}_t(z)=\frac{c_j\boldsymbol{\alpha}_j(t)}{z-z(\gamma_j(t))}+\boldsymbol{\psi}_{j,0}(t)+O\big((z-z(\gamma_j))\big),
\end{equation}
we may see that $\boldsymbol{\psi}_t\mm_t=\zeta_{t,i}(\widehat{p})\cdot\boldsymbol{\psi}_t$ locally. Note that the set $\{\zeta_{t,i}\}$ of local meromorphic functions has poles 
only at the pre-images $\pi^{-1}(nK)$ of the canonical divisor $K$ on $\mathfrak{R}$ and they are invariant under time shift, since $\zeta_{t,i}(\widehat{p})=\widehat{w_i}(\widehat{p})^{m_i}\mu(\widehat{p})^{n_i}$ in the neighborhoods of $\pi^{-1}(p_i)$ where $\pi:\widehat{\mathfrak{R}}\to\mathfrak{R}$ and $\widehat{w_i}$ is the lifting of $w_i$. From Equation~\eqref{coseq1}, we have $(\lambda_t-\zeta_{t,i})\cdot\boldsymbol{\psi}_t=\frac{d}{dt}\boldsymbol{\psi}_t$ around 
$p_i$. What this says is that the poles of $\lambda_t$ at $\pi^{-1}(p_i)$ are also isospectral, since $\frac{d}{dt}\boldsymbol{\psi}_t$ does not have poles at $\pi^{-1}(p_i)$. Consequently, Theorem~\ref{cothe1} confirms the linearity of this flow induced by such $\mm$, since  $\frac{d}{dt}\LL_t=\partial\rho(\mm_t)=\text{ constant}$.

We may see that adding an element in $\HH^0(\widehat{\mathfrak{R}},\pi^\ast\K^n)$ to $\rho(\mm_t)$ is equivalent to adding an element commuting with $\ml$ to $\mm$ in Lax representation~\eqref{de5lax}: 
Consider a time-dependent matrix $\mq_t(p)$ such that $[\mq_t,\ml_t]=0$ where $p\in\mathfrak{R}$. Since $\mq_t$ and $\ml_t$ commute with each other, 
$\mq_t$ preserves the eigenspaces of $\ml_t$. Accordingly, there is a global meromorphic function $\vartheta_t(\hap)$ on $\widehat{\mathfrak{R}}$ such that
\begin{equation}\label{equi10}
\boldsymbol{\psi}_t\mq_t=\vartheta_t(\hap)\cdot\boldsymbol{\psi}_t.
\end{equation}
Moreover, since $\boldsymbol{\psi}_t\ml_t=\mu(\hap)\cdot\boldsymbol{\psi}_t$ and the poles of $\mu(\hap)$ are at $\pi^{-1}(nK)$, we see that the poles of $\vartheta_t$ are only at $\pi^{-1}(nK)$. Thus we conclude that $\vartheta_t\in\HH^0(\widehat{\mathfrak{R}},\pi^\ast\K^n)$. Note that $\vartheta_t$ is not necessarily isospectral unless $\mq_t$ is of form $\mathsf{P}(\ml_t)$ where $\mathsf{P}$ is a polynomial. Combining Equation~\eqref{coseq1} with Equation~\eqref{equi10}, we have
$$\boldsymbol{\psi}_t (\mm_t+\mq_t)+\frac{d}{dt}\boldsymbol{\psi}_t=(\lambda_t+\vartheta_t)\cdot\boldsymbol{\psi}_t.$$
Consequently, we see that $\rho(\mm_t+\mq_t)\equiv\rho(\mm_t)$ modulo $\HH^0(\widehat{\mathfrak{R}},\pi^\ast\K^n)$.

In conclusion, the underlying machinery of this observation is that the sum of residues is zero. More precisely, what this implies is that the behavior of $\widehat{g}+l-1$ poles is translated into the behavior of the lifting divisor in $\pi^{-1}(nK)$. The linearity of the dynamics of $\widehat{g}+l-1$ poles is encoded by the linearity of the dynamics of the lifting divisor in $\pi^{-1}(nK)$. After normalizing by $\mm(p_0)=0$ of form~\eqref{chm}, we denote this straight line flow by $\mathsf{a}=(p_i,n_i,m_i)$. Note that $m_i$ can be a negative integer. It is not hard to see that theses flows commute with each other (\emph{Theorem 2.1} in \cite{kri02}). Moreover, by constructing a symplectic structure on $\mathcal{L}^K/\mathbf{SL}(l,\mathbb{C})$, Krichever calculates Hamiltonians. The Hamiltonian of the flow associated with $\mathsf{a}=(p_i,n_i,m_i)$ is given by

$$H_{\mathsf{a}}(\ml)=-\frac{1}{n_i}\res{p_i}\tr(w^{-m_i}\ml^{n_i})dz\text{ for }\mathsf{a}=(p_i,n_i,m_i)\text{ where }$$
$w_i$ is a local coordinate around $p_i$. See \emph{p.248} in \cite{kri02} for more detailed investigation.

\end{document}